\documentclass[leqno,12pt]{article}

\usepackage{latexsym,amsmath,amsthm,amssymb,amsfonts,amscd}
\usepackage{graphicx} 
\usepackage{pstricks}
\evensidemargin\oddsidemargin

\def \e{\mathbb E}
\def \E {{\rm E}} 
\def \p{\mathbb P}
\renewcommand{\P} {{\rm P}}

\def \U{\mathcal U}
\def \root{ e }
\def \proot{{e_*}}
\def\px{ {x_*}}
\def \pu{ {u_*}} 
\def \tree{{\mathcal T}}
\def \ptree{{\mathcal T}_*}
\def \ray {\mathcal R}
\def\pT{ {T_{*}}}
\def\pU{ \U_{*}}
\def \B{{\mathcal B}}

\def\t{{\mathbb T}}
\def\pt{{\mathbb T}_{*}}

\title{Speed of the biased random walk on a Galton--Watson tree}
\author{Elie A\"id\'ekon\footnote{LPMA, Universit\'e Paris 6, 4, place Jussieu,75005 Paris, France. Email: elie.aidekon@upmc.fr}}

\begin{document}

\baselineskip=18pt
\setcounter{page}{1}

\renewcommand{\theequation}{\thesection.\arabic{equation}}
\newtheorem{theorem}{Theorem}[section]
\newtheorem{lemma}[theorem]{Lemma}
\newtheorem{proposition}[theorem]{Proposition}
\newtheorem{corollary}[theorem]{Corollary}
\newtheorem{remark}[theorem]{Remark}
\newtheorem{fact}[theorem]{Fact}
\newtheorem{problem}[theorem]{Problem}

\maketitle

{\leftskip=2truecm \rightskip=2truecm \baselineskip=15pt \small

\noindent{\slshape\bfseries Summary.} We give an expression of the speed of the biased random walk on a Galton--Watson tree. In the particular case of the simple random walk, we recover the result of Lyons, Pemantle and Peres \cite{LyPePe95}. The proof uses a description of  the invariant distribution of the environment seen from the particle.

\bigskip

\bigskip


} 

\section{Introduction}

Let $\t$ be a Galton--Watson tree with root $\root$, and $\nu$ be its offspring distribution with values in $\mathbb{N}$. We suppose that $m:=\e[\nu]>1$, so that the tree is super-critical. In particular, the event $\mathcal S$ that $\t$ is infinite has a positive probability, and we let $q:=1-\p(\mathcal S)<1$ be the extinction probability. We call $\nu(x)$  the number of children of the vertex $x$ in $\t$. For $x\in \t\backslash\{e\}$,  we denote by $\px$ the parent of $x$, that is the neighbour of $x$ which lies on the path from $x$ to the root $\root$, and by $xi,1\le i\le \nu(x)$ the children of $x$. We call $\pt$ the tree $\t$ on which we add an artificial parent $\proot$ to the root $\root$.

For any $\lambda >0$, and conditionally on $\pt$, we introduce the $\lambda$-biased random walk $(X_n)_{n\ge 0}$ which is the Markov chain such that, for $x\neq \proot$,
\begin{eqnarray}
\P(X_{n+1}=\px\,|\, X_n=x) &=& {\lambda \over \lambda+\nu(x)}, \label{p(x,px)}\\
\P(X_{n+1}=xi\,|\, X_n=x) &=& {1\over \lambda + \nu(x)} \; \; \mbox{for any} \;\; 1\le i\le \nu(x), \label{p(x,xi)}
\end{eqnarray} 
and which is reflected at $\proot$. It is easily seen that this Markov chain is reversible. We denote by $\P_x$ the quenched probability associated to the Markov chain $(X_n)_n$ starting from $x$ and by $\p_x$ the annealed probability obtained by averaging $\P_x$ over the Galton--Watson measure. They are respectively associated to the expectations $\E_x$ and $\e_x$. \\

When $\lambda<m$, we know from Lyons \cite{Lyons} that the walk is almost surely transient on the event $\mathcal S$. Moreover, if we denote by $|x|$ the generation of $x$, Lyons, Pemantle and Peres \cite{LyPePe96} showed that, conditionally on $\mathcal S$, the limit $\ell_\lambda:=\lim_{n\to\infty} {|X_n| \over n}$ exists almost surely, is determinist and is positive if and only if $\lambda \in (\lambda_c,m)$ with $\lambda_c:=\e[\nu q^{\nu-1}]$. This is the regime we are interested in.

For any vertex $x \in \pt$, let 
\begin{equation}\label{def:tau}
\tau_x:=\min\{n\ge 1\,:\, X_n=x\}
\end{equation}

\noindent be the hitting time of the vertex $x$ by the biased random walk, with the notation that $\min\emptyset:=\infty$, and, for $x\neq \proot$,
$$
\beta(x):=\P_x(\tau_{\px}=\infty)
$$
be the quenched probability of never reaching the parent of $x$ when starting from $x$. Notice that we have $\beta(x)>0$ if and only if the subtree rooted at $x$ is infinite. Then, let $(\beta_i,i\ge 0)$ be, under $\p$, generic i.i.d. random variables distributed as $\beta(\root)$, and independent of $\nu$. 
\begin{theorem}\label{t:speed}
Suppose that $m\in(1,\infty)$ and $\lambda\in( \lambda_c  ,m)$. Then,
\begin{equation}\label{eq:l_{lambda}}
\ell_\lambda = \e\left[{(\nu-\lambda) \beta_0 \over \lambda -1+\sum_{i=0}^\nu \beta_i}\right] \Bigg / \e\left[{(\nu+\lambda) \beta_0 \over \lambda -1+\sum_{i=0}^\nu \beta_i}\right].
\end{equation}
\end{theorem}

\noindent Notice that $\ell_{\lambda}$ is the speed of a $\lambda$-biased random walk on a ``regular'' tree with offspring $m_{\lambda}= \e\left[{\nu \beta_0 \over \lambda -1+\sum_{i=0}^\nu \beta_i}\right] / \e\left[{ \beta_0 \over \lambda -1+\sum_{i=0}^\nu \beta_i}\right]$. The FKG inequality implies that $m_{\lambda} \le m$, which means that the randomness of the tree slows down the walk, as conjectured in \cite{LyPePe97}, and already proved in \cite{Chen} and \cite{Virag}.  \\

\noindent The speed in the case $\lambda=1$ was already obtained by Lyons, Pemantle and Peres \cite{LyPePe95}, who found that $\ell_1=\e[{\nu-1\over \nu+1}]$. This can be seen from (\ref{eq:l_{lambda}}) using symmetry. Indeed, taking $\lambda=1$, we see that the numerator is
$\e\left[(\nu-1) { \beta_0 \over \sum_{i=0}^\nu \beta_i}\right] = \e\left[(\nu-1)/(\nu+1)\right]$, while the denominator is just $1$.  In the case $\lambda\to m$, which stands for the near-recurrent regime, Ben Arous, Hu, Olla and Zeitouni \cite{BaHuOlZe} computed the derivative of $\ell_\lambda$, establishing the Einstein relation. Interestingly, the authors give another representation of the speed $\ell_{\lambda}$, at least when $\lambda$ is close enough to $m$. In the zero speed regime $\lambda \le \lambda_c$, Ben Arous, Fribergh, Gantert and Hammond \cite{BaFrGaHa} showed tightness of the properly rescaled random walk, though a limit law fails. A central limit theorem was obtained by Peres and Zeitouni \cite{PeZe06}, by means, in the case $\lambda=m$, of a construction of the invariant distribution on the space of trees. The invariant distribution in the case $\lambda>m$ was given in \cite{BaHuOlZe}. We mention that, so far, the only  case in the transient regime $\lambda<m$ for which such an invariant distribution was known  was the simple random walk case $\lambda=1$ studied in  \cite{LyPePe95}. Theorem \ref{t:inv} in Section \ref{s:inv} gives a description of the invariant measure for all $\lambda\in(\lambda_c ,m)$. These  measures are the limit measures of the tree rooted at the current position of the walker as time goes to infinity. In particular, these measures lie on the space of trees with a backbone, the backbone standing for the  ray linking the walker to the root.  In the setting of random walks on Galton--Watson trees with random conductances, Gantert, M\"uller, Popov and Vachkovskaia  \cite{gantert}  obtained a similar formula for the speed via the construction of the invariant measure in terms of effective conductances. \\

The paper is organized as follows. Section \ref{s:prel} introduces some notation and the concept of backward tree seen from a vertex. Section \ref{s:rever} investigates the law of the tree seen from a vertex that we visit for the first time. Using a time reversal argument, we are able to describe the distribution of this tree in Proposition \ref{p:rever}. Then, we obtain in Section \ref{s:inv} the invariant measure of the tree seen from the particle. Theorem \ref{t:speed} follows in Section \ref{s:speed}.

\section{Preliminaries}
\label{s:prel}

\subsection{The space of words $\U$}
\label{s:words}
We let $\U:=\{\root\}\cup \bigcup_{n\ge 1}(\mathbb{N}^*)^n$ be the set of words, and $|u|$ be the length of the word $u$, where we set  $|\root|:=0$. We equip $\mathcal U$ with the lexicographical order. For any word $u\in \U$ with label $u=i_1\ldots i_n$,  we denote by  $\overline u\in \U$ the word with letters in reversed order $\overline u:= i_n\ldots i_1$ (and $\overline \root:=\root$). If $u\neq \root$, we denote by $\pu$ the parent of $u$, that is the word $i_1\ldots i_{n-1}$, and by $u_{*_k}$ the word $i_1\ldots i_{n-k}$, which stands for the ancestor of $u$ at generation $|u|-k$. We have $u_{*_k}:=\root$ if $k=|u|$ and $u_{*_k}:=u$ if $k=0$. Finally, for $u,v\in\U$,  we denote by $uv$ the concatenation of $u$ and $v$.  We add to the set of words  the element $\proot$, which stands for the parent of the root and we write $\pU:=\U\cup\{\proot\}$. We set $|\proot|=-1$, hence $u_{*_k}=\proot$ for $k=|u|+1$ for any $u\in\U$. We denote by $\ray_x:=\{x_{*_k},\, 1\le k\le |x|+1\}$ the set of strict ancestors of $x$.

\subsection{The space of trees $\tree$}
\label{s:tree}
Following Neveu \cite{Neveu}, a tree $T$ is defined as a subset of $\U$ such that
\begin{itemize}
\item $\root \in T$,
\item if $x\in T\backslash\{\root\}$, then $\px\in T$,
\item if $x=i_1\ldots i_n \in T\backslash\{\root\}$, then any word $i_1\ldots i_{n-1}j$ with $j\le i_n$ belongs to $T$.
\end{itemize}

\noindent We call $\tree$ the space of all trees $T$. For any tree $T$, we define $\pT$ as the tree  on which we add the parent $\proot$ to the root $\root$. Then, let $\ptree:=\{ \pT,T\in\tree\}$. For a tree $T\in\tree$, and a vertex $u\in\pT$, we denote by $\nu_T(u)$ or $\nu_{\pT}(u)$  the number of children of $u$ in $\pT$, and we notice that $\nu_T(\proot)=\nu_{\pT}(\proot)=1$. We will write only $\nu(u)$ when  there is no doubt about which tree we are dealing with. 

\bigskip

We introduce double trees. For any $u\in\U$, let $u^-:=(u,-1)$ and $u^+:=(u,1)$.  Given two trees $T,T^+\in \tree$, we define the double tree $T \!\!- \!\!\!\bullet T^+$ as the tree obtained by drawing an edge between the roots of $T$ and $T^+$. Formally, $T\!\!- \!\!\!\bullet T^+$   is the set $\{ u^-,\,u\in T \} \cup \{ u^+,\, u\in T^+\}$. We root the double tree at $\root^+$. Given $r$ an element of $T$, we say that $X$ is the $r$-parent of $Y$ in $T\!\! - \!\!\!\bullet T^+$ if either
\begin{itemize}
\item  $Y=y^+$  and  $X=y_*^+$,
\item $Y=\root^+$ and $X=\root^-$,
\item $Y=y^-$ with $y\notin \ray_r \cup\{ u\in\U\,:\, u\ge r\}$ and $X=y_*^-$,
\item $Y=r_{*_k}^-$ and $X=r_{*_{k-1}}^-$ for some $k\in[1,|r|]$.
\end{itemize}

\noindent In words, the $r$-parent of a vertex $x$ is the vertex which would be the parent of $x$ if we were ''hanging" the tree at $r$.  Notice that we defined the $r$-parent only for the vertices  which do not belong to $\{ u^-\,:\, u\in\U,\, u\ge r \}$.

\begin{figure}[h]
\begin{center}
\resizebox{9cm}{7cm}{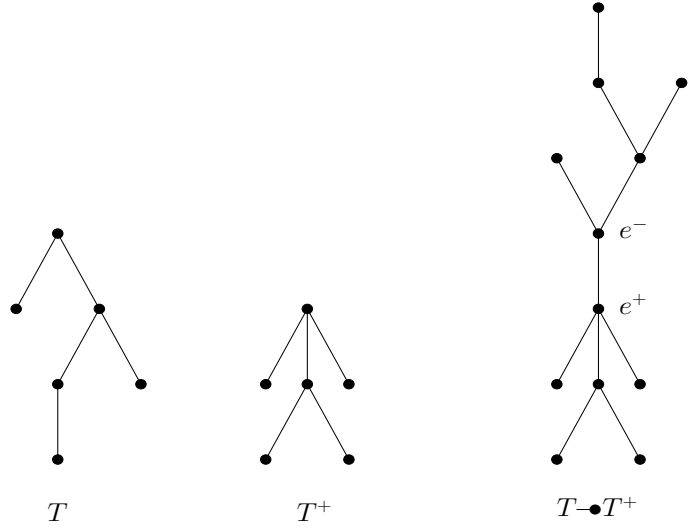}
\end{center}
\caption{A double tree}
\end{figure}

\subsection{The backward tree $\B_x(\pT)$}

Let $\delta$ be some cemetery tree. For a tree $\pT \in \ptree$ and a word $x\in \U$, we define the tree $\pT^{\le x} \in \ptree\cup\{\delta\}$ cut at $x$  by
$$
\pT^{ \le x} := 
\begin{cases}
\delta &\mbox{ if } x\notin \pT,\\
\pT\backslash \{u\in \U\,:\, x<u \} &\mbox{ if } x\in \pT.
\end{cases}
$$

\begin{figure}[h]
\begin{center}
\resizebox{8cm}{7cm}{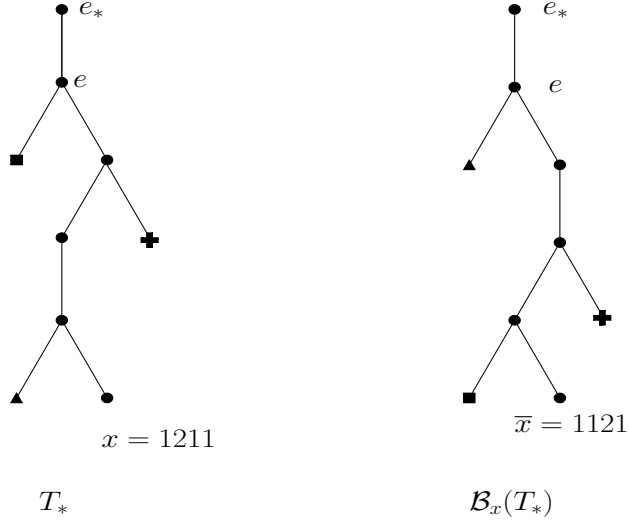}
\end{center}
\caption{The backward tree at $x$}
\end{figure}

\noindent In other words, if $x\in \pT$, then $\pT^{\le x}$ is the tree $\pT$ in which you remove the strict descendants of $x$. We call $\pU^{\le x}$ the set of words $\pU\backslash \{u\in \U\,:\, x<u \}$. We now introduce the backward tree at $x$. For any word $x\in \U$, let $\Psi_x:\pU^{\le x} \to \pU^{\le \overline x}$ such that: 
\begin{itemize}
\item for any $k\in[0,|x|+1]$, $\Psi_x(x_{*_k}) = \overline x_{*_{|x|-k+1}}$,
\item for any $k\in[1,|x|]$ and $v\in \U$ such that $x_{*_k}v$ is not a descendant of $x_{*_{k+1}}$, $\Psi_x(x_{*_k}v)=\Psi_x(x_{*_k})v$. 
\end{itemize}

\noindent The application $\Psi_x$ is a bijection, with inverse map $\Psi_{\overline x}$.  For any tree $\pT\in \ptree$, we call backward tree at $x$ the tree
\begin{equation}\label{def:backward}
\B_x(\pT) := \Psi_x(\pT^{\le x}),
\end{equation}
image of $\pT^{\le x}$ by $\Psi_x$, with the notation that $\Psi_x(\delta):=\delta$. This is the tree obtained by cutting the descendants of $x$ and then ''hanging'' the tree $\pT$ at $x$. We observe that,
\begin{itemize}
\item  $\nu_{\B_x(\pT)}(\proot)=1$,
\item $\nu_{\B_x(\pT)}(\overline x)=0$,
\item for any other $u\in \B_x(\pT)$, we have $\nu_{\B_x(\pT)}(u)=\nu_{\pT}(\Psi_{\overline x}(u))$.
\end{itemize}

\bigskip

\noindent Recall that $\t$ is a Galton--Watson tree with offspring distribution $\nu$. 
\begin{lemma}\label{l:galton}
Let $x\in \U$. The distributions of the trees $\B_x(\pt)$ and $\pt^{\le\overline x}$ are the same.
\end{lemma}

\noindent {\it Proof}. For any sequence $(k_u,u\in\U)\in \mathbb N^{\U}$, denote by  $\mathcal M(k_u,u\in\U) \in \ptree$ the unique tree such that for any $u\in \mathcal M(k_u,u\in\U)$ the number of children of $u$ is $1$ if $u=\proot$ and  $k_u$ otherwise. Take $(\kappa(u),u\in\U)$ i.i.d. random variables distributed as $\nu$. Then notice that the tree $\mathcal M(\kappa(u),u\in\U)$ is distributed as $\pt$. Therefore, we set in this proof 
$$
\pt:=\mathcal M(\kappa(u),u\in\U).
$$ 
We check that we can extend the map $\Psi_{\overline x}$ to a bijection on $\pU$ by letting $\Psi_{\overline x}({\overline x}v):= x v$ for any strict descendant ${\overline x}v$ of $ {\overline x}$.  Suppose that $x\in\pt$.  We know that if $u \in \B_x(\pt)$, then the number of children of $u$ is $1$ if $u=\proot$, $0$ if $u=\overline x$ and $\kappa(\Psi_{\overline x}(u))$ otherwise. By definition, this yields that   
 $$
\B_x(\pt) =\mathcal M(\kappa(\Psi_{\overline x}(u)){\bf 1}_{\{u\neq {\overline x}\}},u\in\U).
 $$
 Let 
$
\widetilde \t_*:= \mathcal M(\kappa(\Psi_{\overline x}(u)),u\in\U)$. We notice that $\mathcal M(\kappa(\Psi_{\overline x}(u)){\bf 1}_{\{u\neq {\overline x}\}},u\in\U)=\widetilde \t_*^{\le \overline x}$.
Therefore, if $x\in\pt$, then
$$
\B_x(\pt)=\widetilde \t_*^{\le \overline x}.
$$

\noindent We check that the equality holds also when $x\notin \pt$. Observe that $\widetilde \t_*$ is distributed as $\pt$ to complete the proof.  \hfill $\Box$

\section{The environment seen from the particle at fresh epochs}
\label{s:rever}

For any tree $\pT \in \ptree$, we denote by $\P^{\pT}$ a probability measure under which $(X_n)_{n\ge 0}$ is a Markov chain on $\pT$  with transition probabilities given by (\ref{p(x,px)}) and (\ref{p(x,xi)}).  For any vertex $x\in \pT$, we denote by $\P^{\pT}_x$ the probability  $\P^{\pT}(\cdot\,|\, X_0=x)$. We will just write $\P_x$ if the tree $\pT$ is clear from the context.

\begin{lemma}\label{l:rever}
Suppose that $\lambda>0$. Let $\pT$ be a tree in $\ptree$, $x$ be a vertex in $\pT\backslash\{\proot\}$ and $(\proot=u_0,u_1,\ldots,u_n=x)$ be a nearest-neighbour trajectory in $\pT$ such that $u_j\notin\{ \proot,x \}$ for any $j\in (0,n)$. Then,
$$
\P^{\pT}_\proot \left( X_j =u_j,\, \forall\, j\le n\right) = \P^{\B_x(\pT)}_\proot\left( X_j =  \Psi_x(u_{n-j})  ,\, \forall\, j\le n\right).
$$
\end{lemma}

\noindent {\it Proof}. We decompose the trajectory $(u_j,j\le n)$ along the ancestral path $\ray_x$. Let $j_0:=0$. Supposing that we know $j_i$, we define $j_{i+1}$ as the smallest integer $j_{i+1}>j_{i}$ such that $u_{j_{i+1}}$ is an ancestor of $x$ different from $u_{j_i}$. Let $m$ be the integer such that $u_{j_{m+1}}=x$. We see that necessarily  $j_1=1$, $(u_{j_0},u_{j_1}) = (\proot,\root)$ and $(u_{j_m},u_{j_{m+1}})=(\px,x)$. For $i\in[1,m]$, let $c_i$ be the cycle $(u_{j_i},u_{j_i+1},\ldots,u_{j_{i+1}-1})$. Notice that in this cycle,   the vertex $u_{j_i}$ is the unique element of $\ray_x$ visited, at least twice at times $j_i$ and $j_{i+1}-1$. We set for any cycle $c=(z_0,z_1,\ldots,z_k)$,
$$
\P^{\pT}(c):=\prod_{\ell=0}^{k-1} \P^{\pT}_{z_\ell}(X_1=z_{\ell+1})
$$

\noindent with the notation that $\prod_{\emptyset}:=1$. Using the Markov property, we see that
\begin{equation}\label{eq:revera}
\P^{\pT}_\proot(X_j =u_j,\, \forall\, j\le n)=\prod_{i=1}^{m} \P^{\pT}(c_i) \prod_{i=1}^m \P^{\pT}_{u_{j_i}}(X_1=u_{j_{i+1}} ).
\end{equation}

\noindent For any vertex $z $, let $a(z):= (\lambda + \nu_\pT(z))^{-1}$. Notice that the term corresponding to $i=m$ in the second product is
$$
\P^{\pT}_\px(X_1= x) = a(\px).
$$

\noindent For any $z\neq \proot$, let $N_u(z)$ be the number of times the oriented edge $(z,z_*)$ is crossed by the trajectory $(u_j,j\le n)$. Notice that the oriented edge $( z_*,z)$ is crossed $1+N_u(z)$ times when $z\in \ray_x$. Using the transition probabilities (\ref{p(x,px)}) and (\ref{p(x,xi)}), we deduce that
$$
 \prod_{i=1}^{m-1} \P^{\pT}_{u_{j_i}}(X_1=u_{j_i+1} )
 =
 \prod_{k=1}^{|x|-1} \left(\lambda a(x_{*_k}) a(x_{*_{k+1}})\right)^{N_u(x_{*_k})}a(x_{*_{k+1}}).
$$

\noindent Therefore, we can rewrite (\ref{eq:revera}) as
\begin{equation}\label{eq:reverb}
\P^{\pT}_\proot(X_j =u_j,\, \forall\, j\le n)=   \mbox{ \Large $\Pi$}_1\mbox{ \Large $\Pi$}_2
\end{equation}
where
\begin{eqnarray} \label{eq:reverb1}
\mbox{ \Large $\Pi$}_1 &:=& \prod_{i=1}^{m} \P^{\pT}(c_i),\\
\mbox{ \Large $\Pi$}_2 &:=& \,a(\px)   \prod_{k=1}^{|x|-1} \left(\lambda a(x_{*_k}) a(x_{*_{k+1}})\right)^{N_u(x_{*_k})}a(x_{*_{k+1}}) . \label{eq:reverb2}
\end{eqnarray}¥

\noindent We look now at the probability $\P^{\B_x(\pT)}_\proot(X_j =v_j,\, \forall\, j\le n)$, where $v_j:=\Psi_x(u_{n-j})$. We decompose the trajectory $(v_j,j\le n)$ along  $\ray_{\overline x}$. Observe that $(v_j,j\le n)$ is the time-reversed trajectory of $(u_j,j\le n)$ looked in the backward tree. Therefore, the cycles of $(v_j,\,j\le n)$ are the image by $\Psi_x$ of the time-reversed cycles of $(u_j,\,j\le n)$. We need some notation. Let ${\buildrel \leftarrow \over c_i}$ be the path $c_i$ time-reversed, and  $\Psi_x\left({\buildrel \leftarrow \over c_i}\right)$ be its image by $\Psi_x$, that is 
$$
\Psi_x\left({\buildrel \leftarrow \over c_i}\right)=(\Psi_x(u_{j_{i+1}-1}),\Psi_x(u_{j_{i+1}-2}), \ldots, \Psi_x(u_{j_i})).$$

\noindent Let 
$$
\P^{\B_x(\pT)}\left(\Psi_x\left({\buildrel \leftarrow \over c_i}\right)\right):=\prod_{\ell=j_i}^{j_{i+1}-2} \P^{\B_x(\pT)}\left(X_1=\Psi_x(u_{\ell}) \,|\, X_0=\Psi_x(u_{\ell+1})\right).
$$

\noindent We introduce for any vertex $z\in\B_x(\pT)$,
$$
a_{\B}(z):= \left(\lambda +\nu_{\B_x(\pT)}(z)\right)^{-1}
$$

\noindent and, for $z\neq \proot$, $N_v(z)$ the number of times the trajectory $(v_j,j\le n)$ crosses the directed edge $(z,z_*)$. Equation (\ref{eq:reverb}) reads for the trajectory $(v_j,\,j\le n)$,
\begin{equation}\label{eq:reverc}
\P^{\B_x(\pT)}_\proot(X_j =v_j,\, \forall\, j\le n) =  \mbox{ \Large $\Pi$}_{\B,1}\mbox{ \Large $\Pi$}_{\B,2} 
\end{equation}

\noindent where
\begin{eqnarray*}
\mbox{ \Large $\Pi$}_{\B,1} &:=& \prod_{i=1}^{m} \P^{\B_x(\pT)}\left(  \Psi_x\left( {\buildrel \leftarrow \over c_i}\right)  \right),\\
\mbox{ \Large $\Pi$}_{\B,2} &:=&  \,a_{\B}({\overline x}_{*}) \prod_{k=1}^{|\overline x|-1} \left(\lambda a_{\B}(\overline x_{*_k}) a_{\B}(\overline x_{*_{k+1}})\right)^{N_v(\overline x_{*_k})}a_{\B}(\overline x_{*_{k+1}}) .
\end{eqnarray*}

\noindent Going from $\pT$ to $\B_x(\pT)$, we did not change the configuration of the subtrees located outside the ancestral path $\ray_x$ of $x$. This yields that $ \P^{\B_x(\pT)}\left( \Psi\left( {\buildrel \leftarrow \over c_i}\right) \right)=\P^{\pT}\left({\buildrel \leftarrow \over c_i}\right)$ which is $ \P^{\pT}(c_i)$ since the Markov chain $(X_n)_{n\ge 0}$ is reversible. By definition of $\mbox{ \Large $\Pi$}_1$ in (\ref{eq:reverb1}), we get
$$
\mbox{ \Large $\Pi$}_{\B,1} = \mbox{ \Large $\Pi$}_1 .
$$

\noindent  We observe that $a_{\B}(z)=a\left(\Psi_{\overline x}(z)\right)$ whenever $z \notin \{\proot,\overline x\}$, and $\Psi_{\overline x}(\overline x_{*_k})= x_{*_{|x|-k+1}}$ by definition. Moreover, for any $k\in[1,|x|-1]$, we have $N_v(\overline x_{*_k})=N_u(x_{*_{|x|-k}})$. This gives that
\begin{eqnarray*}
\mbox{ \Large $\Pi$}_{\B,2} 
&=& 
a(\root) \prod_{k=1}^{|x|-1} \left(\lambda a(x_{*_{| x|-k+1}}) a(x_{*_{| x|-k}})\right)^{N_u(x_{*_{| x|-k}})}a(x_{*_{| x|-k}}) \\
 &=&
a(\px) \prod_{k=1}^{|x|-1} \left(\lambda a(x_{*_{| x|-k+1}}) a(x_{*_{| x|-k}})\right)^{N_u(x_{*_{| x|-k}})}a(x_{*_{| x|-k+1}}) ,
\end{eqnarray*}

\noindent hence, recalling (\ref{eq:reverb2}),  $\mbox{ \Large $\Pi$}_{\B,2} =\mbox{ \Large $\Pi$}_{2} $. Equations (\ref{eq:reverb}) and (\ref{eq:reverc}) lead to
$$
\P^{\pT}_\proot(X_j =u_j,\, \forall\, j\le n) =  \P^{\B_x(\pT)}_\proot(X_j =v_j,\, \forall\, j\le n)
$$
which completes the proof.  \hfill $\Box$

\bigskip

We introduce $\xi_k$, the $k$-th distinct vertex visited by the walk, and $\theta_k:=\tau_{\xi_k}$. These variables are respectively called fresh points, and fresh epochs in \cite{LyPePe96}. They can be defined  by $\theta_0=0$, $\xi_0=X_0$ and for any $k\ge 1$ by
\begin{eqnarray}\label{def:theta}
\theta_k &:=& \min \{ i>\theta_{k-1}  \mbox{ such that } X_i\notin \{X_j,0 \le j<i \}\} , \\
\xi_k &:=& X_{\theta_k}.\label{def:xi}
\end{eqnarray}

\noindent We give the distribution of the tree seen at a fresh epoch $\theta_k$, conditionally  on $\{\theta_k<\tau_{\proot} \}$.

\bigskip

\begin{proposition}\label{p:rever}
Suppose that $\lambda>0$. Let $k\ge 1$. Under $\p_\proot\left(\cdot\,|\, \theta_k<\tau_{\proot}\right)$, we have 
$$
\left(\B_{\xi_k}(\pt),(\Psi_{\xi_k}(X_{\theta_k-j}))_{j\le \theta_k} \right) 
 \, {\buildrel (d) \over =} \, 
 \left(\pt^{\le \xi_k}, (X_j)_{j\le \theta_k}\right).
$$
\end{proposition}

\noindent {\it Proof}. For any relevant bounded measurable map $F$ and any word $x\in\U$, we have
\begin{eqnarray*}
&& \e_{\proot}\left[ F\left( \B_{\xi_k}( \pt), (\Psi_{\xi_k}(X_{\theta_k-j}))_{j\le \theta_k} \right){\bf 1}_{\{ \xi_k=x, \theta_k<\tau_{\proot}\}} \right]\\
&=&
\e_{\proot}\left[ F\left(\B_{x}(\pt),(\Psi_{x}(X_{\theta_k-j}))_{j\le \theta_k} \right){\bf 1}_{\{ \xi_k=x, \theta_k<\tau_{\proot}\}} \right]\\
&=&
\e_{\proot}\left[ F\left(\B_{x}(\pt),(\widetilde X_{j})_{j\le \widetilde \theta_k} \right){\bf 1}_{\{ \widetilde \xi_k=\overline x, \widetilde \theta_k<\widetilde \tau_{\proot}\}} \right]
\end{eqnarray*}
by Lemma \ref{l:rever}, where $(\widetilde X_n)_{n\ge 0}$ is the $\lambda$-biased random walk on the tree $\B_x(\pt)$, and the variables $\widetilde \theta_k$, $\widetilde \xi_k$ and $\widetilde \tau_\proot$ are the analogues of $\theta_k$, $\xi_k$ and $\tau_\proot$ for the Markov chain $(\widetilde X_n)_{n\ge 0}$. By Lemma \ref{l:galton}, it yields that
\begin{eqnarray*}
\e_\proot\left[ F\left(\B_{\xi_k}(\pt),(\Psi_{\xi_k}(X_{\theta_k-j}))_{j\le \theta_k} \right){\bf 1}_{\{ \xi_k=x, \theta_k<\tau_{\proot}\}} \right]
&=&
\e_\proot\left[ F\left( \pt^{\le \overline x},(X_{j})_{j\le \theta_k} \right){\bf 1}_{\{ \xi_k=\overline x, \theta_k<\tau_{\proot}\}} \right]\\
&=&
\e_\proot\left[ F\left(\pt^{\le \xi_k},(X_{j})_{j\le \theta_k} \right){\bf 1}_{\{ \xi_k=\overline x, \theta_k<\tau_{\proot}\}} \right].
\end{eqnarray*}

\noindent We complete the proof by summing over $x\in\U$. \hfill $\Box$

\bigskip 

The last lemma gives the asymptotic probability that $n$ is a fresh epoch. To state it, we introduce the regeneration epochs $(\Gamma_k,k\ge 0)$ defined by $\Gamma_0:=\inf\{ \ell \in \{\theta_k,k\ge 0\}  \,:\,  X_j\neq (X_\ell)_*\, \forall\, j\ge \ell,\,X_\ell\neq \proot   \}$ and for any $k\ge 1$,
\begin{equation}\label{def:gamma}
\Gamma_k := \inf\{ \ell > \Gamma_{k-1}\,:\,  \ell \in \{\theta_k,k\ge 1\}, X_j \neq (X_\ell)_*\, \forall\, j\ge \ell    \},
\end{equation}

\noindent where $(X_\ell)_*$ stands for the parent of the vertex $X_\ell$. For any $k\ge 0$, it is well-known that, under $\p$, the random walk after time $\Gamma_k$ is independent of its past. Moreover, the walk  $(X_\ell,\, \ell\ge \Gamma_k)$ seen in the subtree  rooted at $X_{\Gamma_k}$ is distributed as $(X_\ell,\ell \ge 0)$ under $\p_\root(\cdot\,|\, \tau_{\proot}=\infty)$. We refer to Section 3 of \cite{LyPePe96} for the proof of such facts. We have that $\Gamma_k<\infty$ for any $k\ge 0$ almost surely on the event $\mathcal S$  when $\lambda<m$, and $ \e_\root\left[  \Gamma_1\,|\,\tau_\proot=\infty\right]<\infty$ if and only if $\lambda\in(\lambda_c,m)$.
\begin{lemma}\label{l:firstvisit} Suppose that $m>1$ and $\lambda\in (0,m)$.
We have
$$
\lim_{n\to\infty} \p_\root(n\in \{\theta_k,\, k\ge 0 \},\tau_\proot>n) = {1\over \e_\root\left[  \Gamma_1\,|\,\tau_\proot=\infty\right]}.
$$
\end{lemma}

\noindent {\it Proof}. By the Markov property at time $n$ and the branching property at vertex $X_n$, we observe that
$$
\p_\root(n\in \{\theta_k,\, k\ge 0 \},\tau_\proot>n)\p_\root(\tau_\proot=\infty) = \p_\root(n\in\{ \Gamma_k,k\ge 0\},\, \tau_\proot = \infty)
$$
hence
$$
\p_\root(n\in \{\theta_k,\, k\ge 0 \},\tau_\proot>n)= \p_\root(n\in\{ \Gamma_k,k\ge 0\}\,|\, \tau_\proot = \infty).
$$

\noindent We mention that $\Gamma_0=0$ on the event that $\tau_\proot=\infty$, when starting from $\root$. Since $(\Gamma_{k+1}-\Gamma_k,k\ge 0)$  is a sequence of i.i.d random variables under $\p_\root(\cdot\,|\, \tau_\proot=\infty)$ with mean $\e_\root\left[\Gamma_1\,|\, \tau_\proot=\infty\right]$, the lemma follows from the renewal theorem pp. 360, XI.1 \cite{feller}. \hfill $\Box$


\section{Asymptotic distribution of the environment seen from the particle}
\label{s:inv}

This section is devoted to the asymptotic distribution of the tree seen from the particle. Since $(X_{n})_{n\ge 0}$ is a random walk biased towards the root, it is important to keep track of the root in the  tree seen from $X_n$. Therefore, we will be interested in trees with a marked ray, defined as a couple $(T_{*},R)$ where $T_{*}\in \ptree$, and $R$ is a (finite or infinite)  self-avoiding path of $T_{*}$ starting from the parent of the root $\proot$. We equip the space of trees, resp. the space of marked trees, with the topology generated by finite subtrees, resp. by finite subtrees with a finite ray. They are Polish spaces. \\

For any tree $T\in \tree$ and any $x\in \pT$, let 
$$
T_x:= \{ u\in T\,:\, u \ge x \}.
$$
We recall that we labelled our trees with the space of words $\U$. Remember that $\proot$ has label ${\overline X}_{n}$ in the backward tree $\B_{X_n}(\pt)$.  Recall from Section \ref{s:words} that $\ray_{x}$ stands for the set of words that are strict ancestors of $x$. We are interested in the asymptotic distribution  of $((\B_{X_n}(\pt),\ray_{{\overline X}_{n}}), \t_{X_n})$ in the product topology.  Let $\t$ and $\t^+$ be two independent Galton--Watson trees.   For any tree $\pT\in\ptree$ and any vertex $x\neq \proot$, we can define  $\beta_{\pT}(x)$ as the probability that the biased random walk on $\pT$ never hits $\px$ starting from $x$. We write only $\beta(x)$ when the tree $\pT$ is clear from the context. We write in the following theorem $\nu^+(e):=\nu_{\t^+}(e)$, $\beta(e):=\beta_{\pt}(e)$, $\beta^+(i):=\beta_{\pt^+}(i)$. Finally, conditionally on $\pt$, let $\ray$ be a random ray of $\pt$ with distribution the harmonic measure. It has the law of the almost sure limit of $\ray_{X_{n}}$ as $n\to\infty$, where $(X_{n})_{n\ge 0}$ is the $\lambda$-biased random walk on $\pt$. Observe that $\ray$ is properly defined on the event that $\pt$ is infinite.

\begin{theorem}\label{t:inv}
Suppose that $m\in (1,\infty)$ and $\lambda\in (\lambda_c,m)$. Under $\p_{\proot}(\cdot\,|\,\mathcal S)$, the random variable $((\B_{X_n}(\pt),\ray_{{\overline X}_{n}}),\t_{X_n})$ converges in distribution as $n\to\infty$. The limit distribution has density
\begin{equation} \label{eq:density}
 C_{\lambda}^{-1} {(\lambda +  \nu^+(e)) \beta(e) \over \lambda-1+\beta(e) +\sum_{i=1}^{ \nu^+(e)} \beta^+(i)}
\end{equation}
with respect to $((\pt,\ray),\t^+)$, where $C_{\lambda}$ is the renormalising constant.
\end{theorem}

In the case $\lambda=1$,  the density (\ref{eq:density}) is given by $ C_{1}^{-1} {(1 +  \nu^+(e)) \beta(e) \over \beta(e) +\sum_{i=1}^{ \nu^+(e)} \beta^+(i)}$. If we look at the couple $(\pt,\t^{+})$ as a rooted tree in which the root has $1+\nu_{+}(e)$ children (the tree $\t$ is then a subtree rooted at a vertex of generation $1$), we can take the projection of the invariant measure on the space of unlabeled rooted trees (without marked ray). We recover that the invariant measure is simply the augmented Galton--Watson measure, as proved in \cite{LyPePe95}. This measure is obtained by attaching to the root  $1+\nu$ independent Galton--Watson trees. 

When $\lambda\to m$, the variable $\beta$ converges to $0$. Therefore,  the density (\ref{eq:density}) is equivalent  to $C_{\lambda}^{-1} {m +\nu^+(e) \over m-1} \beta(e)$ as $\lambda\to m$. Proposition 3.1 of \cite{BaHuOlZe} shows that,  when $\nu$ admits a second moment, ${\beta(e) \over \e[\beta]}$ is bounded in $L^2$, which implies that $C_{\lambda}\sim {2m\over m-1}\e[\beta]$, and converges in law. The limit is the distribution of the random variable $W:=\lim_{n\to\infty}  {1\over m^n}\#\{ x\in \t\,:\, |x|=n\}$. Consequently, when $\nu$ has a second moment, the density (\ref{eq:density}) converges in law to ${ m +\nu^+(e) \over 2m} W$ as $\lambda \to m$. This agrees with the invariant measure found in \cite{PeZe06} in the recurrent case $\lambda=m$, and denoted there by {\tt  IGWR}.

\subsection{On the conductance $\beta$}

In this section, let $\pT\in\ptree$ be a fixed tree, and write $\beta(x)$, $\nu(x)$ for  $\beta_{\pT}(x)$, $\nu_{\pT}(x)$. The quantity $\beta(\root)$ is also called conductance of the tree, because of the link between reversible Markov chains and electrical networks, see \cite{DoSn}. It satisfies the recurrence equation 
\begin{equation}\label{eq:beta}
\beta(e) ={\sum_{i=1}^{\nu(e)} \beta(i) \over \lambda + \sum_{i=1}^{\nu(e)} \beta(i) }.
\end{equation}

\noindent  Letting $\beta_n(x)$ be the probability to hit level $n$ before $\px$, we have actually, for $n\ge 1$,
\begin{equation}\label{eq:betan}
\beta_n(e) ={\sum_{i=1}^{\nu(e)} \beta_n(i) \over \lambda + \sum_{i=1}^{\nu(e)} \beta_n(i) }.
\end{equation}

\noindent This is easily seen from the Markov property. Indeed, notice that
\begin{eqnarray*}
\beta_n(e)
=
\sum_{k\ge 0}\P_\root^{\pT}(\tau_\root<\tau_\proot\land \tau_n)^k\P_\root^{\pT}(\tau_n<\tau_\root)
\end{eqnarray*}

\noindent where $\tau_n$ is the hitting time of level $n$. Since 
$$
\P_\root^{\pT}(\tau_\root<\tau_\proot\land \tau_n)=\sum_{i=1}^{\nu(e)}{1\over \lambda + \nu(e)}(1-\beta_n(i))
$$ 

\noindent and
$$
\P_\root^{\pT}(\tau_n<\tau_\root) = \sum_{i=1}^{\nu(e)}{1\over \lambda + \nu(e)}\beta_n(i),
$$ 

\noindent equation (\ref{eq:betan}) follows. Let $n\to\infty$ to get (\ref{eq:beta}). The next lemma implies that the renormalizing constant in Theorem \ref{t:inv} is finite indeed.
\begin{lemma}\label{l:E(beta)}
Suppose that $m>1$ and $\lambda \in (\lambda_c,m)$.  We have
$$
\e\left[ {{\bf 1}_{\mathcal S} \over \lambda - 1 + \beta(e)} \right]<\infty.
$$
\end{lemma}

\noindent {\it Proof}.  The statement is trivial if $\lambda> 1$. Suppose first that $\lambda<1$. By coupling with a one-dimensional random walk, we see that on the event $\mathcal S$, we have $\beta(e)\ge 1-\lambda$. In particular, $\beta_n(e)\ge 1-\lambda$ for any $n\ge 1$.  Use the recurrence equation (\ref{eq:betan}) to get that
\begin{equation}\label{eq:E(beta)1}
 {\beta_n(e) \over \lambda - 1 + \beta_n(e)} = {1\over \lambda}  {\sum_{i=1}^{\nu(e)} \beta_n(i) \over \lambda - 1 + \sum_{i=1}^{\nu(e)} \beta_n(i)}.
\end{equation}

\noindent On the event $\mathcal S$, there exists an index $I\le \nu(e)$ such that the tree rooted at $I$ is infinite. Since $\beta_n(I)\ge 1-\lambda$, we see that
$$
{\sum_{i \le \nu(e),i\neq I} \beta_n(i) \over \lambda - 1 + \sum_{i=1}^{\nu(e)} \beta_n(i)} \le 1.
$$

\noindent On the event that there exists $J\neq I$ such that the tree rooted at $J$ is also infinite, we have
$$
 { \beta_n(I) \over \lambda - 1 + \sum_{i=1}^{\nu(e)} \beta_n(i)} \le {\beta_n(I)\over \beta_n(J)} \le {1\over 1-\lambda}.
$$

\noindent We get that
\begin{eqnarray*}
\e\left[ {\sum_{i=1}^{\nu(e)} \beta_n(i) \over \lambda - 1 + \sum_{i=1}^{\nu(e)} \beta_n(i)}\right]
&\le&
1+{1\over 1-\lambda}
+
\e\left[ {\beta_n(I){\bf 1}_{\{ \beta(j)=0\,\forall j\neq I\}} \over \lambda - 1 + \sum_{i=1}^{\nu(e)} \beta_n(i)}\right]\\
&=&
{\lambda \over 1- \lambda}
+
\e\left[{ \beta_n(I){\bf 1}_{\{ \beta(j)=0 \, \forall j\neq I\}} \over \lambda - 1 +  \beta_n(I) }\right]\\
&=&
{\lambda \over 1- \lambda}
+
\e\left[ \nu q^{\nu-1} \right]\e\left[{ \beta_{n-1}(e)\over \lambda - 1 +  \beta_{n-1}(e) }\right]
.
\end{eqnarray*}

\noindent Recall that $\lambda_c:=\e\left[ \nu q^{\nu-1} \right]$. In view of (\ref{eq:E(beta)1}), we end up with, for any $n\ge 1$,
$$
\e\left[  {\beta_n(e) \over \lambda - 1 + \beta_n(e)}   \right] \le {1\over 1-\lambda} + {\lambda_c\over \lambda}\e\left[ {\beta_{n-1}(e)\over \lambda - 1 +  \beta_{n-1}(e) }\right].
$$

\noindent Applying the above inequality for $n, n-1,\ldots, 1$, we obtain that, for any $\lambda\in(\lambda_c,1)$ and any $n\ge 1$,
$$
\e\left[  {\beta_n(e) \over \lambda - 1 + \beta_n(e)}   \right] \le  {1\over 1-\lambda} {1\over 1- (\lambda_c/\lambda)} + \left({\lambda_c \over \lambda}\right)^n {1\over \lambda}.
$$ 

\noindent Fatou's lemma yields that
$$
\e\left[  {\beta(e) \over \lambda - 1 + \beta(e)}   \right] \le  {1\over 1-\lambda} {1\over 1- (\lambda_c/\lambda)}.
$$ 

\noindent Observe that $\e\left[  {\beta(e) \over \lambda - 1 + \beta(e)}   \right] \ge (1-\lambda)\e\left[  {{\bf 1}_{\mathcal S} \over \lambda - 1 + \beta(e)}   \right]$ to complete the proof in the case $\lambda<1$. In the case $\lambda=1$, we have to show that $\e\left[{{\bf 1}_{\mathcal S}\over \beta(e)}\right]<\infty$. By (\ref{eq:betan}), we have, on the event $\mathcal S$,
$$
{1\over \beta_n(e)} = 1 + {1 \over \sum_{i=1}^{\nu(e)} \beta_n(i)}.
$$

\noindent Let $\varepsilon>0$. With $I$ being defined as before, we check that, on the event $\mathcal S$,
$$
 { 1 \over \sum_{i=1}^{\nu(e)} \beta_n(i)} \le {1\over \beta_n(I)}{\bf 1}_{\{ \beta(i)<\varepsilon\forall\,i\neq I \}} + {1\over \varepsilon}.
$$

\noindent Hence,
$$
\e\left[{{\bf 1}_{\mathcal S}\over \beta_n(e)} \right] \le 1+\varepsilon^{-1} + \e\left[{{\bf 1}_{\mathcal S}\over \beta_{n-1}}\right]\e\left[\nu q_\varepsilon^{\nu-1}\right]
$$

\noindent with $q_\varepsilon:=\p(\beta(e)< \varepsilon)$. Notice that $q_\varepsilon\to q$ as $\varepsilon\to 0$. Taking $\varepsilon>0$ small enough such that $\lambda_\varepsilon:=\e\left[\nu q_\varepsilon^{\nu-1}\right]<1$, we have that
$$
\e\left[{{\bf 1}_{\mathcal S}\over \beta_n(e)} \right] \le(1+\varepsilon^{-1}){1\over 1-\lambda_\varepsilon} + \lambda_{\varepsilon}^{n}(1-q).
$$
Use Fatou's lemma	 to complete the proof.  \hfill $\Box$

\subsection{Random walks on double trees}

\label{s:RWdouble}

Recall that we introduced the concepts of double trees and of $r$-parents in Section \ref{s:tree}. For two trees $T,T^+ \in \tree$, and under some probability $\P^{T\!\! - \!\!\!\bullet T^+}_{\root^+}$,  we introduce two Markov chains on the double tree $T\!\! - \!\!\!\bullet T^+$.\\
For any $r\in T$, we define the biased random walk $(Y_n^{(r)})_{n\ge 0}$ on $T\!\! - \!\!\!\bullet T^+$ with respect to $r$ as the Markov chain, starting from $\root^+$ which moves with weight $\lambda$ to the $r$-parent of the current vertex, with weight $1$ to the other neighbors and which is reflected at the vertex $r^-$. In particular, $Y_n^{(r)}$ never visits the subtree $\{u^-,\, u>r\}$. In words, $(Y_n^{(r)})_{n\ge 0}$ is the $\lambda$-biased random walk on the tree rerooted at $r$.\\
On the other hand, we define $(Y_n)_{n\ge 0}$ the Markov chain on $T\!\! - \!\!\!\bullet T^+$ which has the transition probabilities of the biased random walk in $T$ and in $T^+$. More precisely, if we set $(\proot,-1):=\root^+$ and $(\proot,1):=\root^-$,  the Markov chain $(Y_n)_{n\ge 0}$, while being at $(u,\eta)\in \mathcal U\times\{-1,1\}$,  goes to $(u_*,\eta)$ with weight $\lambda$ and to $(ui,\eta)$ with weight $1$,  this for every child $ui$  of $u$ in $T$ if $\eta=-1$ and every child $ui$ of $u$ in $T^+$ if $\eta=1$. . 

\begin{figure}[ht]
\begin{center}
\resizebox{5cm}{9cm}{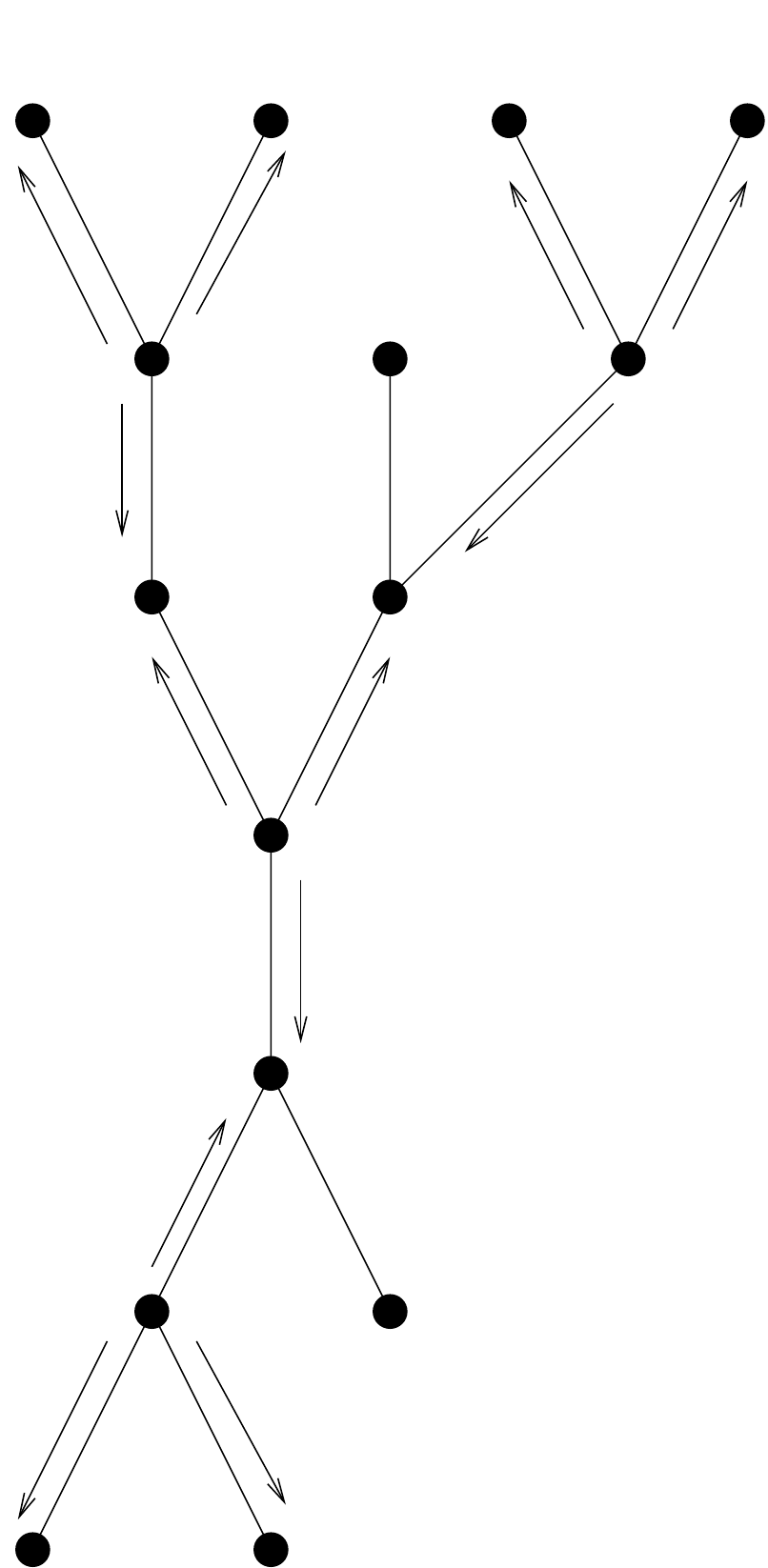}\qquad \qquad \qquad \qquad
\resizebox{5cm}{8,5cm}{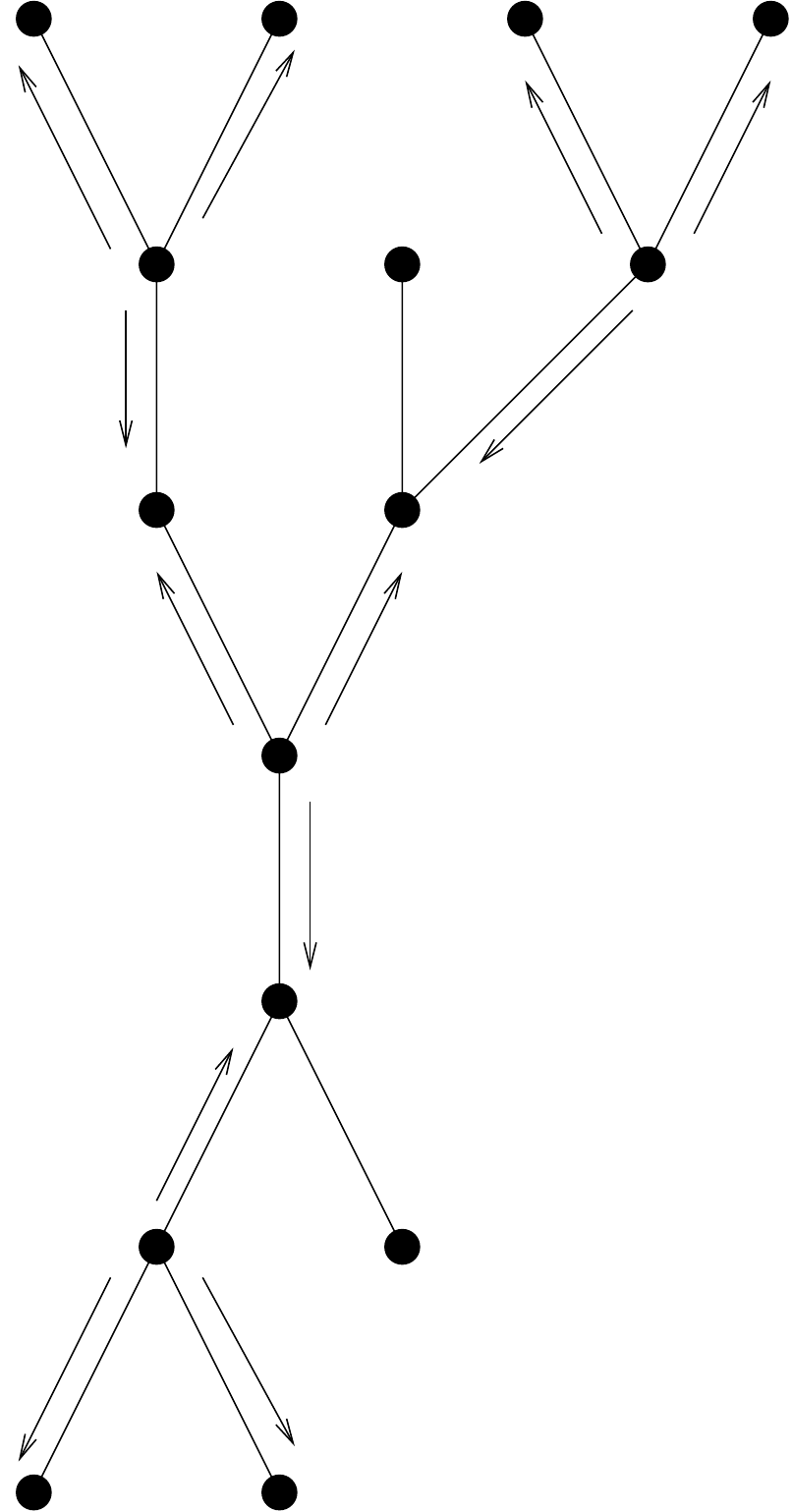}
\end{center}
\caption{The Markov chains  $Y^{(r)}$ (left) and $Y$ (right)}
\end{figure}

\begin{lemma}\label{l:Y(r)}
Let $T\!\! - \!\!\!\bullet T^+$ be a double tree. Let $(\root^+=u_0,u_1,\ldots u_n=\root^+)$ be a sequence of vertices in $T\!\! - \!\!\!\bullet T^+$ such that $u_k\notin \{u^-,\, u\ge r\}$ for any $k\le n$. Denoting by $N_u(y,z)$  the number of crosses of the directed edge $(y,z)$ by the trajectory $(u_k)_{k\le n}$, we have
$$
\P^{T\!\! - \!\!\!\bullet T^+}_{\root^+}\left( Y^{(r)}_k=u_k,\,\forall\,k\le n \right)= \lambda^{-N_u(\root^+,\root^-)}\P^{T\!\! - \!\!\!\bullet T^+}_{\root^+}\left( Y_k=u_k ,\,\forall\,k\le n\right).
$$ 
\end{lemma}

\noindent {\it Proof}.  Let $p^{(r)}(x,y)$, resp. $p(x,y)$, denote the transition probability of the walk $Y^{(r)}$, resp. the walk $Y$, from $x$ to $y$. We have
$$
\P^{T\!\! - \!\!\!\bullet T^+}_{\root^+}\left( Y^{(r)}_k=u_k,\,\forall\,k\le n \right)
=
\prod_{k=0}^{n-1} p^{(r)}(u_k,u_{k+1}).
$$

\noindent Similarly,
$$
\P^{T\!\! - \!\!\!\bullet T^+}_{\root^+}\left( Y_k=u_k,\,\forall\,k\le n \right)
=
\prod_{k=0}^{n-1} p(u_k,u_{k+1}).
$$

\noindent We notice that $p^{(r)}(u_k,u_{k+1})=p(u_k,u_{k+1})$ if $u_k$ or $u_{k+1}$ does not belong to $\{r_{*_\ell}^-,\, \ell \in [1, |r|+1] \}$ where we recall that $\proot^-:=\root^+$. Hence, we only have to show that
\begin{eqnarray}\label{eq:double1a}
&& \prod_{\ell=1}^{|r|} \left(p^{(r)}(r_{*_\ell}^-,r_{*_{\ell+1}}^-)\right)^{N_u(r_{*_{\ell}}^-,r_{*_{\ell+1}}^-)} \left(p^{(r)}(r_{*_{\ell+1}}^-,r_{*_{\ell}}^-)\right)^{N_u(r_{*_{\ell+1}}^-,r_{*_{\ell}}^-)}\\
&=&
\lambda^{-N_u(\root^+,\root^-)} \prod_{\ell=1}^{|r|}      \left(p(r_{*_{\ell}}^-,r_{*_{\ell+1}}^-)\right)^{N_u(r_{*_{\ell}}^-,r_{*_{\ell+1}}^-)} \left(p(r_{*_{\ell+1}}^-,r_{*_{\ell}}^-)\right)^{N_u(r_{*_{\ell+1}}^-,r_{*_{\ell}}^-)} .  \nonumber
\end{eqnarray}

\noindent This comes from the following observations: for any $\ell\in[1,|r|-1]$, $p^{(r)}(r_{*_\ell}^-,r_{*_{\ell+1}}^-) =\lambda^{-1} p(r_{*_{\ell+1}}^-,r_{*_{\ell}}^-)$ and $p^{(r)}(r_{*_{\ell+1}}^-,r_{*_{\ell}}^-)= \lambda p(r_{*_{\ell+1}}^-,r_{*_{\ell}}^-)$.  For $\ell=|r|$, we have $p^{(r)}(r_{*_\ell}^-,r_{*_{\ell+1}}^-) =\lambda^{-1} p(r_{*_{\ell+1}}^-,r_{*_{\ell}}^-)$ and $p^{(r)}(r_{*_{\ell+1}}^-,r_{*_{\ell}}^-)=  p(r_{*_{\ell+1}}^-,r_{*_{\ell}}^-)$. Furthermore, $N_u(r_{*_{\ell}}^-,r_{*_{\ell+1}}^-)=N_u(r_{*_{\ell+1}}^-,r_{*_{\ell}}^-)$ for any $\ell\in [1,|r|]$. A straightforward computation yields (\ref{eq:double1a}), and completes the proof. \hfill $\Box$

\bigskip 

\noindent For any $\ell\ge 0$, let $N_\ell^Y(\root^+,\root^-):=\sum_{k=0}^{\ell-1} {\bf 1}_{\{ Y_k = \root^+,Y_{k+1}=\root^- \}}$ with $\sum_{\emptyset}:=0$. We call $\E^{T\!\! - \!\!\!\bullet T^+}_{\root^+}$ the expectation associated to the probability $\P^{T\!\! - \!\!\!\bullet T^+}_{\root^+}$. In the next lemma, we write $\beta(x)=\beta_{\pT}(x)$, $\beta^+(x)=\beta_{T_*^+}(x)$ and $\nu^+(e)=\nu_{T^+}(e)$.

\begin{lemma}\label{l:Y}
Let $T\!\! - \!\!\!\bullet T^+$ be an infinite double tree. We have
\begin{equation}
\E^{T\!\! - \!\!\!\bullet T^+}_{\root^+}\left[    \sum_{\ell \ge 0}   \lambda^{-N_{\ell}^Y(\root^+,\root^-)} {\bf 1}_{\{ Y_\ell=e^+\}}   \right] ={ \lambda +  \nu^+(e)  \over \lambda-1+\beta(e)+\sum_{i=1}^{\nu^+(e)} \beta^+(i)}.
\end{equation}
\end{lemma}

\noindent {\it Proof}. We compute the left-hand side. We observe that
$$
 \sum_{\ell \ge 0}   \lambda^{-N_{\ell}^Y(\root^+,\root^-)} {\bf 1}_{\{ Y_\ell=e^+\}}
 =
\sum_{k \ge 0}   \lambda^{-k}  \sum_{\ell\ge 0}  {\bf 1}_{\{  N^Y_\ell(\root^+,\root^-) = k, Y_\ell=\root^+ \}}.
$$

\noindent Let $(s_k,\,k\ge 0)$ be the stopping times defined by
$$
s_{k}:=\inf\{ \ell  \ge 0\,:\,   N_\ell^Y(\root^+,\root^-)=k \}.
$$ 

\noindent We define  $t_k:=\inf\{\ell \ge s_k\,:\, X_\ell =\root^+  \}$, and we have that $t_0=s_0=0$. Notice that, for any $k\ge 0$,
$$
\sum_{\ell\ge 0}  {\bf 1}_{\{  N^Y_\ell(\root^+,\root^-) = k, Y_\ell=\root^+ \}} 
=
{\bf 1}_{\{t_k<\infty\}} \sum_{\ell=t_k}^{s_{k+1}}  {\bf 1}_{\{  Y_\ell=\root^+ \}}.
$$

\noindent This gives that
$$
\E^{T\!\! - \!\!\!\bullet T^+}_{\root^+}\left[ \sum_{\ell \ge 0}   \lambda^{-N_{\ell}^Y(\root^+,\root^-)} {\bf 1}_{\{ Y_\ell=e^+\}} \right]
 =
\sum_{k \ge 0}   \lambda^{-k} \E^{T\!\! - \!\!\!\bullet T^+}_{\root^+}\left[{\bf 1}_{\{t_k<\infty\}} \sum_{\ell=t_k}^{s_{k+1}}  {\bf 1}_{\{  Y_\ell=\root^+ \}}\right].
$$

\noindent By the strong Markov property at time $t_k$, we have, for any $k\ge 0$,
$$
 \E^{T\!\! - \!\!\!\bullet T^+}_{\root^+}\left[{\bf 1}_{\{t_k<\infty\}}\sum_{\ell=t_k}^{s_{k+1}}  {\bf 1}_{\{  Y_\ell=\root^+ \}}\right]
 =
  \P^{T\!\! - \!\!\!\bullet T^+}_{\root^+}(t_k<\infty)\E^{T\!\! - \!\!\!\bullet T^+}_{\root^+}\left[\sum_{\ell=0}^{s_{1}}  {\bf 1}_{\{  Y_\ell=\root^+ \}}\right].
$$

\noindent We see that $\P^{T\!\! - \!\!\!\bullet T^+}_{\root^+}(t_k<\infty) = \left[(1-\beta^+(\root))(1-\beta(\root))\right]^k$. Moreover, for $\tau_{\root^+}^Y:=\inf\{n\ge 1\,:\, Y_n=\root^+\}$, we have $\P^{T\!\! - \!\!\!\bullet T^+}_{\root^+}(\tau^Y_{\root^+} < s_1)  =  {1\over \lambda + \nu^+(\root)}\sum_{i=1}^{\nu^+(e)} (1-\beta^+(i))$.  This yields that
$$
\E^{T\!\! - \!\!\!\bullet T^+}_{\root^+}\left[\sum_{\ell=0}^{s_{1}}  {\bf 1}_{\{  Y_\ell=\root^+ \}}\right] = {1\over 1-\P^{T\!\! - \!\!\!\bullet T^+}_{\root^+}(\tau^Y_{\root^+}<s_1)}
=
{\lambda +\nu^+(e) \over \lambda + \sum_{i=1}^{\nu^+(\root)} \beta^+(i)}.
$$

\noindent Since $T\!\! - \!\!\!\bullet T^+$ is infinite, we have by coupling with a one-dimensional random walk, $\beta(e)> 1-\lambda$ or $\beta^+(e)> 1-\lambda$. Hence $\lambda^{-1}(1-\beta^+(e))(1-\beta(e))< 1$. We end up with
$$
\E^{T\!\! - \!\!\!\bullet T^+}_{\root^+}\left[ \sum_{\ell \ge 0}   \lambda^{-N_{\ell}^Y(\root^+,\root^-)} {\bf 1}_{\{ Y_\ell=e^+\}} \right]
=
{1\over 1- \lambda^{-1}(1-\beta(\root))(1-\beta^+(\root))}{\lambda + \nu^+(\root) \over \lambda + \sum_{i=1}^{\nu^+(\root)} \beta^+(i)}.
$$

\noindent Apply the recurrence equation (\ref{eq:beta}) to $\beta^+(e)$ to complete the proof. \hfill $\Box$

\bigskip 

\subsection{Proof of Theorem \ref{t:inv}}

{\it Proof of Theorem \ref{t:inv}}. Let $F_1$ and $F_2$ be two bounded measurable functions respectively on the space of marked trees and on $\tree$ which depend only on a finite subtree. Recall the definition of the regeneration epochs $(\Gamma_k,k\ge 1)$ in (\ref{def:gamma}). We will show that
 \begin{eqnarray}\label{eq:invmain}
&& \lim_{n\to\infty}\e_\proot \left[  F_1\left( \B_{X_n}(\pt),\ray_{{\overline X}_{n}}\right)  F_2\left(\t_{X_n}\right) {\bf 1}_{\mathcal S}\right] \\
\nonumber &=&
    {\p(\mathcal S)\over \e_e[\Gamma_1{\bf 1}_{\{ \tau_\proot=\infty\}}]}  \e \left[  F_1( \pt,\ray)F_2(\t^+)  {(\lambda + \nu^+(e) )\beta(e) \over \lambda-1+\beta(e)+\sum_{i=1}^{\nu^+(e)} \beta^+(i)} \right]
\end{eqnarray}

\noindent which proves the theorem. Let us prove (\ref{eq:invmain}).   We first show that 
 \begin{eqnarray}\label{eq:inv1}
&&\lim_{n\to\infty}\e_\proot \left[  F_1\left( \B_{X_n}(\pt),\ray_{{\overline X}_{n}}\right)  F_2\left(\t_{X_n}\right)   {\bf 1}_{\{\tau_{\proot}>n\}} \right] \\
&=&
\nonumber    {1\over \e_e[\Gamma_1\,|\, \tau_\proot=\infty]}  \e \left[  F_1( \pt,\ray)F_2(\t^+)   {(\lambda + \nu^+(e))\beta(e)  \over \lambda-1+\beta(e)+\sum_{i=1}^{\nu^+(e)} \beta^+(i)} \right]. 
\end{eqnarray}

\noindent Let $\varepsilon\in (0,1)$ and, for any random tree $T$, $\mathcal S_{T}$ be the event that $T$ is infinite.  We deduce from dominated convergence that
\begin{eqnarray}\label{eq:inv1a}
&&\e_\proot \left[  F_1\left( \B_{X_n}(\pt),\ray_{{\overline X}_{n}}\right)  F_2\left(\t_{X_n}\right)   {\bf 1}_{\{\tau_{\proot}>n\}} \right] \\
&=&
\e_\proot \left[  F_1\left( \B_{X_n}(\pt),\ray_{{\overline X}_{n}}\right)  F_2\left(\t_{X_n}\right)   {\bf 1}_{\{\tau_{\proot}>n ,|X_n|\ge n^\varepsilon \}}{\bf 1}_{\mathcal S_{\B_{X_n}(\pt)}} \right] +o_n(1). \nonumber 
\end{eqnarray}

\noindent    Recall  the definition of $\theta_k$ and $\xi_k$ in (\ref{def:theta}) and (\ref{def:xi}). We have for any $n\ge 1$,
\begin{eqnarray*}
&& \e_\proot \left[  F_1\left( \B_{X_n}(\pt),\ray_{{\overline X}_{n}}\right)  F_2\left(\t_{X_n}\right)   {\bf 1}_{\{\tau_{\proot}>n , |X_n|\ge n^{\varepsilon} \}} {\bf 1}_{\mathcal S_{\B_{X_n}(\pt)}}\right] \\
&=&
\sum_{k\ge 1} \e_\proot\left[F_1\left( \B_{\xi_k}(\pt),\ray_{{\overline \xi}_{k}} \right)   F_2\left( \t_{\xi_k} \right)  {\bf 1}_{\{X_n=\xi_k,\tau_{\proot}>n,|\xi_k|\ge n^\varepsilon\}}{\bf 1}_{\mathcal S_{\B_{X_n}(\pt)}} \right].
\end{eqnarray*}

\noindent  We want to reroot the tree at $\xi_k$. Notice that $\t_{\xi_k}$ is a Galton--Watson tree independent of $\B_{\xi_k}(\pt)$. By the strong Markov property at time $\theta_k$ and  Proposition \ref{p:rever}, we have that for any $k\ge 1$,
\begin{eqnarray*}
&& \e_\proot\left[F_1\left( \B_{\xi_k}(\pt),\ray_{{\overline \xi}_{k}} \right)   F_2\left( \t_{\xi_k} \right)  {\bf 1}_{\{X_n=\xi_k,\tau_{\proot}>n,|\xi_k|\ge n^\varepsilon\}} {\bf 1}_{\mathcal S_{\B_{X_n}(\pt)}} \right] \\
&=&
\e_{\proot}\left[  F_1(\pt^{\le \xi_k},\ray_{\xi_{k}})F_2(\t^+)    {\bf 1}_{\{ Y_{n-\theta_k}^{(\xi_k)}=e^+,\tau_{\xi_k}^{(\xi_k)}>n-\theta_k \}} {\bf 1}_{\{\tau_{\proot}>\theta_k , |\xi_k|\ge n^\varepsilon \}}   {\bf 1}_{\mathcal S_{\pt^{\le \xi_k}}}\right].
\end{eqnarray*}

\noindent In the last expectation, the Markov chain $(X_n)_{n\ge 0}$ being the biased random walk on $\pt$ starting at $\proot$, the variables $\theta_k$, $\xi_k$ and $\tau_x$ are given by (\ref{def:theta}), (\ref{def:xi}) and (\ref{def:tau}). Moreover, conditionally on $\t$, $\t^+$ and $\{X_\ell,\ell\le \theta_k\}$, we take $(Y_n^{(\xi_k)})_{n\ge 0}$ a biased random walk starting at $\root^+$ with respect to $\xi_k$ on the double tree ${\t\!\! - \!\!\!\bullet \t^+}$ as defined in Section \ref{s:RWdouble}, and $\tau_{\xi_k}^{(\xi_k)} := \inf\{ \ell\ge 1\,:\, Y_\ell^{(\xi_k)}=(\xi_k,-1)\}$. Since $F_1$ depends only on a finite subtree, we get that for $n$ large enough,
\begin{eqnarray}\label{eq:inv1b}
&& \e_\proot \left[  F_1\left( \B_{X_n}(\pt) ,\ray_{{\overline X}_{n}}\right)  F_2\left(\t_{X_n}\right)   {\bf 1}_{\{\tau_{\proot}>n , |X_n|\ge n^{\varepsilon} \}}{\bf 1}_{\mathcal S_{\B_{X_n}(\pt)}}  \right] \\
&=&
\sum_{k\ge 1}\e_{\proot}\left[  F_1(\pt,\ray_{\xi_{k}})F_2(\t^+)    {\bf 1}_{\{ Y_{n-\theta_k}^{(\xi_k)}=e^+,\tau_{\xi_k}^{(\xi_k)}>n-\theta_k \}} {\bf 1}_{\{\tau_{\proot}>\theta_k , |\xi_k|\ge n^\varepsilon \}}  {\bf 1}_{\mathcal S_{\pt^{\le \xi_k}}} \right].\nonumber
\end{eqnarray}

\noindent Lemma \ref{l:Y(r)} implies that
\begin{eqnarray}\label{eq:inv1c}
&&
\sum_{k\ge 1}\e_{\proot}\left[  F_1(\pt,\ray_{\xi_{k}})F_2(\t^+)    {\bf 1}_{\{ Y_{n-\theta_k}^{(\xi_k)}=e^+,\tau_{\xi_k}^{(\xi_k)}>n-\theta_k \}} {\bf 1}_{\{\tau_{\proot}>\theta_k , |\xi_k|\ge n^\varepsilon \}}  \right]\\
\nonumber &=&
\sum_{k\ge 1}\e_{\proot}\left[F_1( \pt,\ray_{\xi_{k}} )F_2(\t^+)     \lambda^{-N_{n-\theta_k}^Y(\root^+,\root^-)} {\bf 1}_{\{ Y_{n-\theta_k}=e^+,\tau_{\xi_k}^Y>n-\theta_k \}}  {\bf 1}_{\{\tau_{\proot}>\theta_k,|\xi_k|\ge n^\varepsilon\}}  {\bf 1}_{\mathcal S_{\pt^{\le \xi_k}}}  \right]
\end{eqnarray}

\noindent where, conditionally on $\t$, $\t^+$, the Markov chain $(Y_n)_{n\ge 0}$ is the biased random walk  on the double tree ${\t\!\! - \!\!\!\bullet \t^+}$ as defined in Section \ref{s:RWdouble}, taken independent of $(X_n)_{n\ge 0}$,  and $\tau_{\xi_k}^{Y} := \inf\{ \ell\ge 1\,:\, Y_\ell=(\xi_k,-1)\}$.
\noindent In view of (\ref{eq:inv1a}), (\ref{eq:inv1b}) and (\ref{eq:inv1c}), we see that, as $n\to\infty$,
\begin{eqnarray*}
&&\e_\proot \left[  F_1\left( \B_{X_n}(\pt),\ray_{{\overline X}_{n}}\right)  F_2\left(\t_{X_n}\right)   {\bf 1}_{\{\tau_{\proot}>n\}} \right] \\
&=&
\e_{\proot}\left[F_2(\t^+)   \sum_{k\ge 1}  F_1( \pt ,\ray_{\xi_{k}}) \lambda^{-N_{n-\theta_k}^Y(\root^+,\root^-)} {\bf 1}_{\{ Y_{n-\theta_k}=e^+,\tau_{\xi_k}^Y>n-\theta_k \}}  {\bf 1}_{\{\tau_{\proot}>\theta_k,|\xi_k|\ge n^\varepsilon\}} {\bf 1}_{\mathcal S_{\pt^{\le \xi_k}}}  \right]\\ 
&& + o_n(1).
\end{eqnarray*}

\noindent Reasoning on the value of $n-\theta_k$, and since $\xi_k=X_{\theta_k}$, we observe that
\begin{eqnarray*}
&&   \sum_{k\ge 1}  F_1( \pt ,\ray_{\xi_{k}})  \lambda^{-N_{n-\theta_k}^Y(\root^+,\root^-)} {\bf 1}_{\{ Y_{n-\theta_k}=e^+,\tau_{\xi_k}^Y>n-\theta_k \}}  {\bf 1}_{\{\tau_{\proot}>\theta_k,|\xi_k|\ge n^\varepsilon\}}  {\bf 1}_{\mathcal S_{\pt^{\le \xi_k}}}   \\
&=&    \sum_{\ell = 0}^{n-1}   F_1( \pt ,\ray_{X_{n-\ell}}) \lambda^{-N_{\ell}^Y(\root^+,\root^-)} {\bf 1}_{\{ Y_\ell=e^+\}} {\bf 1}_{\{\tau_{\proot}>n-\ell,|X_{n-\ell}|\ge n^\varepsilon,\tau_{X_{n-\ell}}^Y>\ell,n-\ell\in\{\theta_k,k\ge 1\}\}} {\bf 1}_{\mathcal S_{\pt^{\le X_{n-\ell}}}} .
\end{eqnarray*}

\noindent Lemma \ref{l:E(beta)} shows that
$$
 \e \left[    {(\lambda + \nu^+(e)){\bf 1}_{\mathcal S_{\t}} \over \lambda-1+\beta(e)+\sum_{i=1}^{\nu^+(e)} \beta^+(i)} \right]<\infty.
$$

\noindent Together with Lemma \ref{l:Y}, it implies that
$$
\e\left[    \sum_{\ell \ge 0}   \lambda^{-N_{\ell}^Y(\root^+,\root^-)} {\bf 1}_{\{ Y_\ell=e^+\}} {\bf 1}_{\mathcal S_{\t}}\right]<\infty.
$$

\noindent  Therefore, we can use dominated convergence to replace 
$$
\sum_{\ell = 0}^{n-1}   F_1( \pt ,\ray_{X_{n-\ell}}) \lambda^{-N_{\ell}^Y(\root^+,\root^-)} {\bf 1}_{\{ Y_\ell=e^+\}} {\bf 1}_{\{\tau_{\proot}>n-\ell,|X_{n-\ell}|\ge n^\varepsilon,\tau_{X_{n-\ell}}^Y>\ell,n-\ell\in\{\theta_k,k\ge 1\}\}} {\bf 1}_{\mathcal S_{\pt^{\le X_{n-\ell}}}}
$$

\noindent by
$$
\sum_{\ell \ge 0}   F_1( \pt ,\ray) \lambda^{-N_{\ell}^Y(\root^+,\root^-)} {\bf 1}_{\{ Y_\ell=e^+\}} {\bf 1}_{\{\tau_{\proot}=\infty,n-\ell\in\{\theta_k,k\ge 1\}\}}
$$

\noindent and hence see that 
\begin{eqnarray*}
&&\e_\proot \left[  F_1\left( \B_{X_n}(\pt),\ray_{{\overline X}_{n}}\right)  F_2\left(\t_{X_n}\right)   {\bf 1}_{\{\tau_{\proot}>n\}} \right] \\
&=&
\e_{\proot}\left[   F_1( \pt ,\ray)  F_2(\t^+)   \sum_{\ell \ge 0}   \lambda^{-N_{\ell}^Y(\root^+,\root^-)} {\bf 1}_{\{ Y_\ell=e^+\}} {\bf 1}_{\{\tau_{\proot}=\infty,n-\ell\in\{\theta_k,k\ge 1\}\}}    \right] + o_n(1).
\end{eqnarray*}

\noindent We deduce from dominated convergence that for any integer $K\ge 1$, we have as well
\begin{eqnarray*}
&&  \e_\proot \left[  F_1\left( \B_{X_n}(\pt),\ray_{{\overline X}_{n}}\right)  F_2\left(\t_{X_n}\right)   {\bf 1}_{\{\tau_{\proot}>n\}} \right] \\
&=&
 \e_{\proot}   \left[  F_1( \pt,\ray )F_2(\t^+)   \sum_{\ell \ge 0}   \lambda^{-N_{\ell}^Y(\root^+,\root^-)} {\bf 1}_{\{ Y_\ell=e^+\}} {\bf 1}_{\{\tau_{\proot}=\infty,n-\ell\in\{\theta_k,k\ge 1\} ,n-\ell\ge \Gamma_K \}}  \right] + o_n(1).
\end{eqnarray*}

\noindent We choose $K$  a deterministic integer such that $F_1$ does not depend on the set $\{u\in \U\,:\, |u|\ge K-1\}$. Notice that necessarily, $|X_{\Gamma_K}|\ge K-1$. In particular, $F_1(\pt,\ray)$  is independent of the subtree rooted at $X_{\Gamma_K}$. Recall that $\t^+$ is independent of $\pt$, hence of $(X_n)_n$ as well.  Using the regenerative structure of the walk $(X_n)_n$ at time $\Gamma_K$, we get that 
\begin{eqnarray*}
&&
 \e_{\proot}   \left[  F_1( \pt,\ray )F_2(\t^+)   \sum_{\ell \ge 0}   \lambda^{-N_{\ell}^Y(\root^+,\root^-)} {\bf 1}_{\{ Y_\ell=e^+\}} {\bf 1}_{\{\tau_{\proot}=\infty,n-\ell\in\{\theta_k,k\ge 1\} ,n-\ell\ge \Gamma_K \}}  \right] \\
 &=&
  \e_{\proot}   \left[  F_1( \pt,\ray )F_2(\t^+)   \sum_{\ell \ge 0}   \lambda^{-N_{\ell}^Y(\root^+,\root^-)} {\bf 1}_{\{ Y_\ell=e^+\}} {\bf 1}_{\{\tau_{\proot}=\infty,n-\ell\ge \Gamma_K\}}  b_{n-\ell-\Gamma_K}\right]
  \end{eqnarray*}
  
\noindent with, for any integer $i \ge 0$,
$  b_i
 : =
  \p_\root\left( i\in \{\theta_k,k\ge 0\}  \, | \, \tau_\proot=\infty  \right)
$.  Lemma \ref{l:firstvisit} says that $b_i\to {1\over \e_e[\Gamma_1\,|\, \tau_\proot=\infty]}   $ as $i\to\infty$, hence
\begin{eqnarray*}
&&
 \lim_{n\to\infty} \e_{\proot}   \left[  F_1( \pt,\ray )F_2(\t^+)   \sum_{\ell \ge 0}   \lambda^{-N_{\ell}^Y(\root^+,\root^-)} {\bf 1}_{\{ Y_\ell=e^+\}} {\bf 1}_{\{\tau_{\proot}=\infty,n-\ell\in\{\theta_k,k\ge 1\} ,n-\ell\ge \Gamma_K \}}  \right] \\
 &=&  {1\over \e_e[\Gamma_1\,|\, \tau_\proot=\infty]}   \e_{\proot}   \left[  F_1( \pt,\ray )F_2(\t^+)   \sum_{\ell \ge 0}   \lambda^{-N_{\ell}^Y(\root^+,\root^-)} {\bf 1}_{\{ Y_\ell=e^+\}} {\bf 1}_{\{\tau_{\proot}=\infty\}} \right].
  \end{eqnarray*}

\noindent Consequently,
\begin{eqnarray*}
&& \lim_{n\to\infty} \e_\proot \left[  F_1\left( \B_{X_n}(\pt),\ray_{{\overline X}_{n}}\right)  F_2\left(\t_{X_n}\right)   {\bf 1}_{\{\tau_{\proot}>n\}} \right]\\
  &=&  {1\over \e_e[\Gamma_1\,|\, \tau_\proot=\infty]}   \e_{\proot}   \left[  F_1( \pt,\ray )F_2(\t^+)   \sum_{\ell \ge 0}   \lambda^{-N_{\ell}^Y(\root^+,\root^-)} {\bf 1}_{\{ Y_\ell=e^+\}} {\bf 1}_{\{\tau_{\proot}=\infty\}} \right].
  \end{eqnarray*}

\noindent Recall that $\beta(e)=\P^{\pt}_\proot(\tau_\proot=\infty)$ by definition. Then apply Lemma \ref{l:Y} to complete the proof of (\ref{eq:inv1}). It remains to remove the conditioning on $\{\tau_{\proot}>n\}$ on the left-hand side. Fix $\ell\ge 1$. For $n\ge\ell$, we have by the Markov property,
$$
\e_\proot\left[F_1(    \B_{X_n}(\pt),\ray_{{\overline X}_{n}} )F_2(  \t_{X_n} )    {\bf 1}_{\{\Gamma_0=\ell\}}\right]  =  \e_\proot \left[ {\bf 1}_{E_\ell} \phi(X_\ell,n-\ell)\right]
$$

\noindent where, for any $k\ge 0$ and $x\in\pt$,
$$
\phi(x,k):=\E_{x} \left[F_1( \B_{X_{k}}(\pt) ,\ray_{{\overline X}_{k}}  )F_2(\t_{X_{k}}  )     {\bf 1}_{\{\tau_{\proot}=\infty\}}\right]
$$ 

\noindent and, for any $\ell\ge 0$, $E_\ell$ is the event that $X_\ell\neq \proot$ and that at time $\ell$, every (non-directed) edge that has been visited at least twice, except the edge between $X_\ell$ and its parent. Since $F_1$ depends on a finite subtree, we can use, when $|X_{n-\ell}|$ is big enough (actually greater than $K-1$), the branching property for the Galton--Watson tree at the vertex  $X_\ell$ to obtain that
$$
\e_\proot\left[ {\bf 1}_{E_\ell} \phi(X_\ell,n-\ell)\right]
=
\p_\proot\left(E_\ell\right) \e_{\root} \left[F_1( \B_{X_{n-\ell}}(\pt),\ray_{{\overline X}_{n-\ell}}   )F_2(\t_{X_{n-\ell}}  )     {\bf 1}_{\{\tau_{\proot}=\infty\}}\right]     + o_n(1).
$$

\noindent Notice that, for any $n-\ell\ge 0$,
\begin{eqnarray*}
&&\e_{\root} \left[F_1( \B_{X_{n-\ell}}(\pt) ,\ray_{{\overline X}_{n-\ell}}  )F_2(\t_{X_{n-\ell}}  )     {\bf 1}_{\{\tau_{\proot}=\infty\}}\right] \\ 
&=&
\e_{\proot} \left[F_1( \B_{X_{n-\ell+1}}(\pt) ,\ray_{{\overline X}_{n-\ell+1}}  )F_2(\t_{X_{n-\ell+1}}  )     {\bf 1}_{\{\tau_{\proot}=\infty\}}\right]. 
\end{eqnarray*}

\noindent Equation (\ref{eq:inv1}) implies that
\begin{eqnarray*}
&& \lim_{n\to\infty} \e_\proot\left[F_1(    \B_{X_n}(\pt),\ray_{{\overline X}_{n}} )F_2(  \t_{X_n} )    {\bf 1}_{\{\Gamma_0=\ell\}}\right] \\
&=&
 \p_\proot(E_\ell) {1\over \e_e[\Gamma_1\,|\, \tau_\proot=\infty]}  \e\left[  F_1( \pt,\ray)F_2(\t^+)   {(\lambda + \nu^+(e))\beta(e)  \over \lambda-1+\beta(e)+\sum_{i=1}^{\nu^+(e)} \beta^+(i)} \right]. 
\end{eqnarray*}

\noindent Since $\{\Gamma_0<\infty\}=\mathcal S$, we deduce that
\begin{eqnarray*}
&& \lim_{n\to\infty} \e_\root\left[F_1(    \B_{X_n}(\pt) ,\ray_{{\overline X}_{n}})F_2(  \t_{X_n} ) {\bf 1}_{\mathcal S}\right] \\
&=&
  { \sum_{\ell\ge 1} \p_\proot(E_\ell) \over \e_e[\Gamma_1\,|\, \tau_\proot=\infty]}  \e_\proot \left[  F_1( \pt,\ray)F_2(\t^+)   {\lambda + \nu^+(e)  \over \lambda-1+\beta(e)+\sum_{i=1}^{\nu^+(e)} \beta^+(i)} \right]. 
\end{eqnarray*}

\noindent We notice that $\p_\proot(E_\ell)\p_\root(\tau_\proot=\infty)=\p_\proot(\Gamma_0=\ell)$, hence
$$
\sum_{\ell\ge 1} \p_e(E_\ell) ={\p(\mathcal S)\over \p_\root(\tau_\proot=\infty)}.
$$

\noindent This proves (\ref{eq:invmain}), hence the theorem.
\hfill $\Box$

\section{Proof of Theorem \ref{t:speed}}
\label{s:speed}

\noindent {\it Proof}.  By dominated convergence, we have $\ell_\lambda = \lim_{n\to\infty} \e_\proot\left[ {|X_n|\over n}\,|\, \mathcal S \right]$. We observe that
$$
 \e_\proot\left[ |X_n| \,|\, \mathcal S  \right] =  \sum_{k=0}^{n-1} \e_\proot\left[  |X_{k+1}| -|X_k| \,|\,\mathcal S \right]= \sum_{k=0}^{n-1} \e_\proot\left[  {\nu(X_k) - \lambda \over \nu(X_k) + \lambda}  \,|\, \mathcal S\right] .
$$

\noindent Use Theorem \ref{t:inv} to complete the proof. \hfill $\Box$

\end{document}